\documentclass[leqno]{article}
\title{
 \bf Rigidity of automorphisms and spherical CR structures}
\author{
  Jih-Hsin Cheng
\thanks{
Research supported in part by National Science Council grant NSC 
87-2115-M-001-006 (R.O.C.).}\\
Institute of Mathematics, Academia Sinica, Taipei}
\date{}

\begin{document}
\maketitle
\begin{abstract}
We establish Bochner-type formulas for operators related
to $CR$ automorphisms and spherical $CR$ structures. From
such formulas, we draw conclusions about rigidity by making
assumptions on the Tanaka-Webster curvature and torsion. 
\end{abstract}

%\setcounter{section}{-1}

%\subjclass
%\endsubjclass
%\keywords
%\endkeywords
%\thanks
%\endthanks
%\abstract
%\endabstract

\section{Statement of results}
It is now clear (e.g. [CL1], [CL2], [CT], [Rum]) that certain distinguished
second-order partial differential operators and their fourth-order
``Laplacians'' play
important roles in the study of three-dimensional $CR$ geometry. In
this paper we will establish Bochner-type formulas for operators 
related to $CR$ automorphisms and spherical $CR$ structures. From
such formulas, we can draw conclusions about rigidity by making
assumptions on the so-called Tanaka-Webster curvature and torsion.

To be precise, let $(M,J,\theta )$ be a smooth, closed (compact without
boundary) 3-dimensional strictly pseudoconvex pseudohermitian manifold.
(see, e.g. [Web],[Le1]) Here $J$ denotes a $CR$ structure and $\theta$
is a contact form, i.e. a nonvanishing real 1-form defining the underlying
contact structure. Associated to $(M,J,\theta )$, we have a canonical affine
connection and notions of Tanaka-Webster scalar curvature and torsion,
denoted $R$ and $A$ (a tensor with coefficient $A_{11}$ or ${A_1}^{\bar 1}$)
respectively. ([Le1],[Tan],[Web]) 

Let $T$ be the characteristic vector field of $\theta$ defined by
${\theta}(T)=1, {\mathcal{L}}_{T}{\theta}=0$. (some authors call $T$ the 
Reeb vector field) By choosing suitable complex vector fields $Z_{1},Z_{\bar 1}
$ such that $JZ_{1}=iZ_{1}, JZ_{\bar 1}=-iZ_{\bar 1}$, we form a (unitary)
frame $\{ Z_{1},Z_{\bar 1},T \}$. ($h_{1{\bar 1}}=1$) The covariant derivatives
are taken with respect to this frame and indicated by ${_{,1}},
{_{,{\bar 1}}},{_{,0}}$ and so on. Let $Aut_{0}(J)$ denote the identity
component of the $CR$ automorphism group with respect to $J$.

Now we can state our first result.

\medskip

$\bf {Theorem \: A}$. Let $(M,J,\theta )$ be a smooth, connected, closed
3-dimensional strictly pseudoconvex pseudohermitian manifold.

(a) Suppose $R < 0, {\sqrt 3}R_{,0}-2Im(A_{11,{\bar 1}{\bar 1}}) > 0$.

\noindent Then $Aut_{0}(J)$ consists of only the identity diffeomorphism.

(b) Suppose $R < 0, {\sqrt 3}R_{,0}-2Im(A_{11,{\bar 1}{\bar 1}}) = 0$

\noindent Then $dim Aut_{0}(J){\leq}1$.

\medskip

We remark that the torsion $A=0$ implies the second condition in (b) holds
due to the Bianchi identity: $R_{,0}=A_{11,{\bar 1}{\bar 1}}+
A_{{\bar 1}{\bar 1},11}$. And the result is compatible with Prop.4.8(a) in
[CT]. On the other hand, we do not know any examples satisfying curvature
conditions in (a).

The idea of proving $Theorem \: A$ goes as follows. Consider a certain
second-order linear partial differential operator ${\mathcal D}_J$ and its ``Laplacian'' ${\mathcal D}_J
^{\star}{\mathcal D}_J$ acting on functions. Here ${\mathcal D}_J^{\star}$
denotes the adjoint of ${\mathcal D}_J$.  The kernel of ${\mathcal D}_J$
parametrizes infinitesimal $CR$ automorphisms. (see [CL1] where ${\mathcal D}_J$ was denoted $B_J^{'}$) Establish a suitable Bochner-type formula for
$<{\mathcal D}_J^{\star}{\mathcal D}_{J}f,f>$ by using commutation relations
and integration by parts repeatedly. (see section 2 for details) $Theorem\: A$
then follows easily from the final formula. Note that the kernel of
${\mathcal D}_J^{\star}$ parametrizes the infinitesimal slice in the study
of $CR$ moduli spaces. ([CL2])

Next we consider deformation of spherical $CR$ structures. A $CR$ structure
or $CR$ manifold is called spherical if it is locally $CR$-equivalent to  
the unit sphere with the standard $CR$ structure. (e.g. [BS],[CL1])
In dimension 3 it can be characterized quantitatively by the vanishing of
a certain fourth-order partial differential operator, the so-called Cartan
(curvature) tensor, denoted $Q_J$. (while in higher dimensions, Chern's
curvature tensor plays the similar role ([CM]), which is of second order)

A spherical $CR$ structure $J$ is called rigid if there is no infinitesimal
deformation up to diffeomorphisms, i.e. for any smooth family of spherical
$CR$ structures $J_{(t)}$ on the base manifold $M$ with $J_{(0)}=J$, 
$d/dt|_{t=0}J_{(t)}$ equals ${\mathcal L}_{X}J$ for some vector field $X$ of
$M$.

To study the rigidity of spherical $CR$ structures, we consider
the linearization of $Q_J$ plus a symmetry-breaking term provided by
${\mathcal D}_{J}{\mathcal D}_J^{\star}$. In section 3 we work out a 
Bochner-type formula and analyze it to obtain pointwise conditions for
$J$ to be rigid.

\medskip

$\bf {Theorem \: B}$. Let $(M,J)$ be a smooth, closed, spherical $CR$ 3-manifold.
Suppose there is a contact form $\theta$ such that $R > 0$, (3.11) and (3.12)
hold. Then $J$ is rigid. 

\medskip 

Note that conditions (a),(b) in $Theorem \: A$ and (3.11),(3.12) in $Theorem \: B$ are independent of positive constant multiples of $\theta$. When the torsion
vanishes, we have the simplified expression as follows.

\medskip

$Corollary \: C$. Let $(M,J)$ be a smooth, closed, spherical $CR$ 3-manifold.
Suppose there is a contact form $\theta$ such that the torsion $A=0$ and

$$
R > 0, 4R(5R^{2}+3{\bigtriangleup}_{b}R)-3|{\nabla}_{b}R|_{\theta}^{2}>0
$$

\noindent Then $J$ is rigid.

\medskip

The sublaplacian and subgradient operators ${\bigtriangleup}_{b},{\nabla}_{b}$
acting on (smooth) functions are defined by
${\bigtriangleup}_{b}f=-{f_{,1}}^{1}
-{f_{,{\bar 1}}}^{\bar 1}$ and ${\nabla}_{b}f={f_{,}}^{1}Z_{1}+{f_{,}}^{\bar 1}
Z_{\bar 1}$ respectively. (cf. [Le1] or [Che]) Also we define $|{\nabla}_{b}f|_{\theta}^{2}=2f_{,1}{f_{,}}^{1}$ for real f.

Observe that $A=0$ and $R$ being
a positive constant satisfy conditions in $Corollary \, C$. In this case the
universal cover of $(M,J)$ must be compact by Rumin's pseudohermitian version
of Myers' theorem ([Rum]), and hence $CR$-equivalent to the standard $S^{3}$. It follows that the fundamental group $\Gamma$ of
$M$ is finite. Hence the group cohomology $H^{1}({\Gamma},{\mathcal G})$ in
deformation theory (p.232 in [BS]) vanishes. So in this special case our
result is compatible with the result obtained by ``Lie theoretical'' argument.
Note that a small perturbation of $A=0$ and $R$ being a positive constant
still satisfies the conditions in $Theorem \: B$.

\medskip

$\bf{Acknowledgment}$. This work was carried out during the author's 
visit at Harvard University. He would therefore like to thank the members
of the Mathematics Department, and especially Professor Shing-Tung Yau,
for their hospitality during his stay.
  
\bigskip

\section{Proof of Theorem A}

Let $\{ {\theta}^{1},{\theta}^{\bar 1},{\theta}\}$ be the coframe dual to
the ``unitary'' frame $\{ Z_{1},Z_{\bar 1},T \}$. (with $h^{1{\bar 1}}=
h_{1{\bar 1}}=1$ in mind, hereafter, we'll write tensors with only lower
indices) Recall ([CL2] or [CL1]) that ${\mathcal{D}}_{J}f= 2Re[(f_{,11}+
iA_{11}f){\theta}^{1}{\otimes}Z_{\bar 1}]$ and the adjoint operator  
${\mathcal D}_J^{\star}$ acts on a deformation tensor $E= 2Re(E_{11}{\theta}^
{1}{\otimes}Z_{\bar 1})$ by

\begin{equation}
{\mathcal D}_J^{\star}E=E_{11,{\bar 1}{\bar 1}}+iA_{11}E_{{\bar 1}{\bar 1}}+
conjugate.
\end{equation}

Also the generalized Folland-Stein operator $L_{\alpha}$ is defined by
$L_{\alpha}f={\bigtriangleup}_{b}f+i{\alpha}f_{,0}$ for a function f.
By Lemma 2.1 in [CL2], we have

\begin{equation}
{\mathcal D}_J^{\star}{\mathcal D}_{J}=(1/2)L_{\alpha}^{\star}L_{\alpha}
+{\mathcal O}_{2}
\end{equation}

\noindent with ${\alpha}=i{\sqrt 3}$, where ${\mathcal O}_{2}$
is an operator of 
weight ${\leq}2$. In the rest of this section, we'll look into the details
of ${\mathcal O}_{2}$. (Note that $L_{\alpha}$ in the leading term of (2.2)
is subelliptic. This implies the existence of the ``infinitesimal slice
decompositions'' in the study of $CR$ moduli spaces without knowing
details of the lower-weight term ([CL2]))

A direct computation shows that for a real-valued function f, we have

\begin{eqnarray}
&&{\mathcal D}_J^{\star}{\mathcal D}_{J}f=f_{,11{\bar 1}{\bar 1}}+f_{,{\bar 1}{\bar 1}11}
+2Re[2iA_{11}f_{,{\bar 1}{\bar 1}} \\
&&+2iA_{11,{\bar 1}}f_{,{\bar 1}}+
iA_{11,{\bar 1}{\bar 1}}f+|A_{11}|^{2}f]. \nonumber
\end{eqnarray}

We'll frequently use the following commutation relations.

\medskip

$\bf {Lemma \: 2.1}$(Ricci identities in pseudohermitian geometry) Let $c_I$
be a coefficient of some tensor with multi-indices $I$. Suppose $I$ consists
of only $1$ and $\bar 1$, and ${\alpha}$=(\# of $1$ in $I$)-(\# of $\bar 1$ in
$I$).Then

\begin{equation}
c_{I,1{\bar 1}}-c_{I,{\bar 1}1}=ic_{I,0}+{\alpha}c_{I}R
\end{equation}

\begin{equation}
c_{I,01}-c_{I,10}=c_{I,{\bar 1}}A_{11}-{\alpha}c_{I}A_{11,{\bar 1}}
\end{equation}

\begin{equation}
c_{I,0{\bar 1}}-c_{I,{\bar 1}0}=c_{I,1}A_{{\bar 1}{\bar 1}}
+{\alpha}c_{I}A_{{\bar 1}{\bar 1},1}
\end{equation}

\noindent
(this lemma generalizes Lemma 2.3 in [Le2] for the three-dimensional case)

\medskip

By using (2.4),(2.5),(2.6) repeatedly, we obtain

\begin{equation}
{\frac{1}{2}}L_{\alpha}^{\star}L_{\alpha}f=f_{,11{\bar 1}{\bar 1}}
+f_{,{\bar 1}{\bar 1}11}
+2Re[({\sqrt 3}-i)(A_{{\bar 1}{\bar 1}}f_{,1})_{,1}-(Rf_{,1})_{,{\bar 1}}]
\end{equation}

\noindent with the choice of ${\alpha}=i{\sqrt 3}$ eliminating
terms having covariant
derivative in the direction $T$. Comparing (2.3) with (2.7) and taking $L^{2}$
-inner product with $f$ give

\begin{eqnarray}
&{\|}{\mathcal D}_{J}f{\|}^{2}&={\frac{1}{2}}{\|}L_{i{\sqrt 3}}f{\|}^{2}
-{\int_M}R|{\nabla}_{b}f|_{\theta}^{2}dv_{\theta} \\
&& +2{\int_M}[Re({\frac{{\sqrt 3}+i}{2}}A_{11,{\bar 1}{\bar 1}})+|A_{11}|^{2}]
f^{2}dv_{\theta} \nonumber \\ 
&& +2{\int_M}Re[(-i-{\sqrt 3})A_{{\bar 1}{\bar 1}}f_{,11}]f
dv_{\theta} \nonumber
\end{eqnarray}

Here the volume form $dv_{\theta}={\theta}{\wedge}d{\theta}$. With respect to
$dv_{\theta}$, we have the divergence theorem and hence integration by parts
in calculus of pseudohermitian geometry. ([Le2],[Che]) For instance,

\begin{eqnarray*}
&{\int_M}A_{{\bar 1}{\bar 1},1}f_{,1}{\bar f}dv_{\theta} & ={\frac{1}{2}}
{\int_M}A_{{\bar 1}{\bar 1},1}(f^{2})_{,1}dv_{\theta}{\:}(f{\:}being{\:}real)\\
 & & =-{\frac{1}{2}}{\int_M}A_{{\bar 1}{\bar 1},11}f^{2}dv_{\theta}
{\:}(integration{\:}by{\:}parts)
\end{eqnarray*}

\noindent was used in deducing (2.8). Now suppose ${\phi}_t{\in}Aut_{0}(J)$ is
a smooth family of $CR$ automorphisms with ${\phi}_{0}=identity$. Then
$X={\frac{d}{dt}}|_{t=0}{\phi}_{t}{\in}$ Lie $Aut_{0}(J)$ (=Lie algebra of
$Aut_{0}(J)$) is an infinitesimal $CR$ automorphism, in particular, an
infinitesimal contact automorphism. According to Lemma 3.4 and 3.5 in 
[CL1], $X=X_f$ is determined by a function $f=-{\theta}(X)$ and satisfies
the relation: ${\mathcal L}_{X}J=2{\mathcal D}_{J}f$. Since ${\mathcal L}_{X}J=0$, we
get

\begin{equation}
0={\mathcal D}_{J}f = f_{,11}+iA_{11}f,
\end{equation}

\noindent and hence 

\begin{equation}
A_{{\bar 1}{\bar 1}}f_{,11} = -i|A_{11}|^{2}f.
\end{equation}

Substituting (2.9), (2.10) in (2.8), we finally obtain

\begin{eqnarray}
& 0 & ={\|}L_{i{\sqrt 3}}f{\|}^{2}-2{\int_M}R|{\nabla}_{b}f|_{\theta}^{2}
dv_{\theta}  \\
& & +{\int_M}[{\sqrt 3}R_{,0}+i(A_{11,{\bar 1}{\bar 1}}-A_{{\bar 1}{\bar 1}
,11})]f^{2}dv_{\theta} \nonumber
\end{eqnarray}

Now it is easy to see from (2.11) that the condition in (a) of $Theorem{\:}A$
implies $f=0$. Therefore $X=X_{f}=0$. For (b), the condition implies
${\nabla}_{b}f=0$. So $f_{,0}=0$ by (2.4). Thus $f$ is constant since $M$
is connected. It follows that $dim(Aut_{0}(J))=dim(Lie{\,}Aut_{0}(J))
{\leq}1$.

\bigskip

\section{Proof of Theorem B}
\setcounter{equation}{0}

 Recall ([CL1]) that the Cartan tensor $Q_{J}=iQ_{11}{\theta}^{1}{\otimes}
Z_{\bar 1}-iQ_{{\bar 1}{\bar 1}}{\theta}^{\bar 1}{\otimes}Z_{1}$ where

$$
Q_{11}={\frac{1}{6}}R_{,11}+\frac{i}{2}RA_{11}-A_{11,0}-\frac{2i}{3}A_{11,
{\bar 1}1}
$$

\noindent and $J$ is spherical if and only if $Q_{J}=0$. (note that we have lowered
indices using $h_{1{\bar 1}}=1$; also $Q_J$ changes ``tensorially'' when
we make a different choice of contact form) Let ${\tilde J}_{(t)}$ be
a smooth family of spherical $CR$ structures with ${\tilde J}_{(0)}=J$.
By a theorem of Gray ([Gra] or [Ham]), there exists a smooth family of
diffeomorphisms $\phi_t$ with $\phi_{0}=identity$ so that for all $t$,
$J_{(t)}={\phi_t^{\star}}{\tilde J}_{(t)}$ has the same underlying
contact structure as $J$ does. Write infinitesimal deformation

$$
\frac{d}{dt}|_{t=0}J_{(t)}= 2E= 4Re(E_{11}{\theta^1}{\otimes}Z_{{\bar 1}})
$$

\noindent and compute $DQ_{J}(2E)={\partial_t}Q_{J_{(t)}}|_{t=0}$ as we did in [CL1].
There appears a ``bad'' term $E_{{\bar 1}{\bar 1},1111}$ in the formula,
so we add a ``symmetry-breaking'' term ${\mathcal D}_{J}{\mathcal D}_J^{\star}
E$ to cancel it. The final formula including terms of lower weights reads

\begin{eqnarray}
& &-DQ_{J}(2E)+\frac{1}{6}{\mathcal D}_{J}{\mathcal D}_J^{\star}
E  =  2Re\{ \frac{1}{3}E_{11,{\bar 1}{\bar 1}11}\\
& &-E_{11,00}
-\frac{2i}{3}E_{11,0{\bar 1}1}
  +  \frac{i}{3}(A_{11}E_{{\bar 1}{\bar 1}})_{,11} \nonumber \\
& &-\frac{1}{6}E_{11}
R_{,1{\bar 1}}+\frac{1}{6}E_{11,{\bar 1}}R_{,1}-\frac{1}{6}
(E_{11}R_{,{\bar 1}})_{,1} \nonumber \\
 & & + \frac{1}{2} A_{11}(iE_{11,{\bar 1}{\bar 1}}-iE_{{\bar 1}{\bar 1},11}
-A_{11}E_{{\bar 1}{\bar 1}}-A_{{\bar 1}{\bar 1}}E_{11})+\frac{i}{2}RE_{11,0}
\nonumber \\
 & & +2A_{11}(A_{11}E_{{\bar 1}{\bar 1}}+A_{{\bar 1}{\bar 1}}E_{11})+
\frac{2i}{3}E_{11}A_{11,{\bar 1}{\bar 1}}-\frac{2i}{3}E_{11,{\bar 1}}
A_{11,{\bar 1}}-\frac{2i}{3}(E_{{\bar 1}{\bar 1}}A_{11,1})_{,1}
\nonumber \\
 & & -\frac{4i}{3}(E_{{\bar 1}{\bar 1},1}A_{11})_{,1}+\frac{i}{6}A_{11}
(E_{11,{\bar 1}{\bar 1}}+E_{{\bar 1}{\bar 1},11}+iA_{11}E_{{\bar 1}{\bar 1}}
-iA_{{\bar 1}{\bar 1}}E_{11})\} {\theta^1}{\otimes}Z_{\bar 1} \nonumber
\end{eqnarray}

The right-hand side of (3.1) can be written as $\frac{1}{12}L_{\alpha}^{\star}
L_{\alpha}E+{\mathcal O}_{2}(E)$ with ${\mathcal O}_{2}$ being an operator
of weight $\leq 2$ and ${\alpha}=4+i{\sqrt 3}$. Since $L_{\alpha}$ is
subelliptic, the above expression was used in [CL1] to show the short-
time solution of a certain regularized evolution equation. Using Lemma 2.1
repeatedly, we can write the highest-weight term of (3.1) as follows:

\begin{eqnarray}
&& \frac{1}{3}E_{11,{\bar 1}{\bar 1}11}-E_{11,00}-\frac{2i}{3}E_{11,0{\bar 1}1}=\frac{1}{3}E_{11,{\bar 1}1{\bar 1}1}-E_{11,00}-iE_{11,0{\bar 1}1}\\
&& +\frac{i}{3}(E_{11,1}A_{{\bar 1}{\bar 1}})_{,1}+\frac{2i}{3}
(E_{11}A_{{\bar 1}{\bar 1},1})_{,1}-\frac{1}{3}(RE_{11,{\bar 1}})_{,1}
\nonumber
\end{eqnarray}
 
On the other hand, we compute

\begin{eqnarray}
&&{\int_M}R|E_{11,1}|^{2}dv_{\theta} = -{\int_M}(RE_{11,1{\bar 1}}
E_{{\bar 1}{\bar 1}}+R_{,{\bar 1}}E_{11,1}E_{{\bar 1}{\bar 1}})dv_{\theta}
\\
&&(by{\:}integration{\:}by{\:}parts)  \nonumber \\
&& =-{\int_M}[R(E_{11,{\bar 1}1}+iE_{11,0}+2RE_{11})E_{{\bar 1}{\bar 1}}
+R_{,{\bar 1}}E_{11,1}E_{{\bar 1}{\bar 1}}]dv_{\theta} 
{\:}(by{\:}(2.4)) \nonumber
\end{eqnarray}

To see how we treat (3.1) in general, we first deal with the torsion$=0$
case. By the local slice theorem ([CL2]), there exists a smooth family of
contact diffeomorphisms $\psi_t$ with ${\psi_0}=identity$ so that 
$J_{(t)}^{'}={\psi}_{t}^{\star}J_{(t)}$ lies in the local slice passing
through $J$. Since the infinitesimal slice at $J$ is parametrized by the
kernel of ${\mathcal D}_J^{\star}$, we have ${\mathcal D}_J^{\star}E=0$
for $2E=\frac{d}{dt}|_{t=0}J_{(t)}^{'}$.

Now applying (3.1) to such a deformation tensor $E$: ${\mathcal D}_J^{\star}
E=0$, substituting (3.2) in (3.1), and taking $L^2$-inner product with $E$,
we obtain

\begin{eqnarray}
& 0 & =<-DQ_{J}(2E),E>{\:}(since{\:}Q_{J_{(t)}^{'}}=0) \\
& & ={\int_M}2Re\{ \frac{1}{3}|E_{11,{\bar 1}1}|^{2}+|E_{11,0}|^{2}
-iE_{11,0}E_{{\bar 1}{\bar 1},1{\bar 1}}+\frac{1}{6}R|E_{11,1}|^{2}
\nonumber \\
&& +\frac{2i}{3}RE_{11,0}E_{{\bar 1}{\bar 1}} 
 -\frac{1}{6}RE_{11,{\bar 1}1}E_{{\bar 1}{\bar 1}}+\frac{1}{6}
R_{,{\bar 1}}E_{11,1}E_{{\bar 1}{\bar 1}} \nonumber \\
&& +[\frac{11}{6}R_{,{\bar 1}1}-
2R_{,1{\bar 1}}+\frac{1}{3}R^{2}]|E_{11}|^{2}\} dv_{\theta} 
\nonumber 
\end{eqnarray}

\noindent by putting $A_{11}=0$, using (3.3) and integration by parts
repeatedly.

Note that we reduce (3.4) to the formula (6.3) in [CL1] for $R$ being a
constant $\hat R$. Furthermore, if ${\hat R}>0$, the right-hand side of
(3.4) is a positive definite quadratic hermitian form in $E_{11,{\bar 1}1}$,
$iE_{11,0}$, and ${\hat R}E_{11}$.

It follows that $E=0$ and $0=\frac{d}{dt}|_{t=0}J_{(t)}^{'}=
\frac{d}{dt}|_{t=0}({\phi_t}{\circ}{\psi_t})^
{\star}{\tilde J}_{(t)}={\mathcal L}_{X}J+\frac{d}{dt}|_{t=0}{\tilde J}_{(t)}$
where the vector field $X=\frac{d}{dt}|_{t=0}({\phi_t}{\circ}{\psi_t})$. 
So $J$ is rigid. For general $R$, we require that the integrand in (3.4)
is a (pointwise) positive definite quadratic hermitian form in 
$E_{11,{\bar 1}1}$, $iE_{11,0}$, $E_{11,1}$, $E_{11}$.
Now the conditions in $Corollary{\:}C$ can be deduced from basic linear
algebra.

When the torsion does not vanish, the formula for $<-DQ_{J}(2E),E>$ with
${\mathcal D}_J^{\star}E=0$ reads

\begin{eqnarray}
& & 0=<-DQ_{J}(2E),E>={\int_M}\{ \frac{2}{3}|E_{11,{\bar 1}1}|^{2}+
2|E_{11,0}|^{2} \\
& & +\frac{1}{3}R|E_{11,1}|^{2} 
+ [\frac{2}{3}R^{2}+\frac{1}{6}{\bigtriangleup}_{b}R+6|A_{11}|^{2}
+\frac{8i}{3}(A_{11,{\bar 1}{\bar 1}}-A_{{\bar 1}{\bar 1},11})]|E_{11}|^{2}
\nonumber \\
& & +2Re[-iE_{11,0}E_{{\bar 1}{\bar 1},1{\bar 1}} 
+\frac{2i}{3}RE_{11,0}E_{{\bar 1}{\bar 1}}-\frac{1}{6}RE_{11,{\bar 1}1}
E_{{\bar 1}{\bar 1}} \nonumber \\
& & +(\frac{1}{6}R_{,{\bar 1}}-2iA_{{\bar 1}{\bar 1},1})
E_{11,1}E_{{\bar 1}{\bar 1}}-\frac{2i}{3}A_{11,1}E_{{\bar 1}{\bar 1},1}
E_{{\bar 1}{\bar 1}}-\frac{5i}{3}A_{11}E_{11,{\bar 1}}
E_{{\bar 1}{\bar 1},{\bar 1}}]\} dv_{\theta} \nonumber
\end{eqnarray}

There are non-cross terms like $A_{11,1}E_{{\bar 1}{\bar 1},1}E_{{\bar 1}{\bar 1}}$ and $A_{11}E_{11,{\bar 1}}E_{{\bar 1}{\bar 1},{\bar 1}}$ in (3.5).
We need an inequality to deal with  $A_{11}E_{11,{\bar 1}}E_{{\bar 1}{\bar 1},{\bar 1}}$.

\medskip

$\bf{Lemma \: 3.1}$. Let $\lambda$, $\rho$ be real numbers. Then

\begin{eqnarray}
2{\lambda}{\rho}{\int_M}R|E_{11,{\bar 1}}|^{2}dv_{\theta} & \leq &
{\lambda}^{2}{\int_M}|E_{11,{\bar 1}1}|^{2}dv_{\theta}+
{\rho}^{2}{\int_M}R^{2}|E_{11}|^{2}dv_{\theta} \\
 & & -{\lambda}{\rho}{\int_M}(R_{,1}E_{11,{\bar 1}}E_{{\bar 1}{\bar 1}}+
R_{,{\bar 1}}E_{{\bar 1}{\bar 1},1}E_{11})dv_{\theta} 
\nonumber
\end{eqnarray}

For the term $A_{11,1}E_{{\bar 1}{\bar 1},1}E_{{\bar 1}{\bar 1}}$, we use
the following estimate:

\begin{equation}
2Re(-\frac{2i}{3}A_{11,1}E_{{\bar 1}{\bar 1},1}E_{{\bar 1}{\bar 1}})
{\geq}-\frac{2}{3}|A_{11,1}|^{\frac{2}{3}}|E_{{\bar 1}{\bar 1},1}|^{2}
-\frac{2}{3}|A_{11,1}|^{\frac{4}{3}}|E_{11}|^{2}.
\end{equation}

\noindent (To deduce (3.7), replace $a$, $b$ by ${\omega}a$, ${\omega}^{2}b$ in the
basic inequality $2Re(ab){\geq}-|a|^{2}-|b|^{2}$ with
$a=-iE_{{\bar 1}{\bar 1},1}$, $b=E_{{\bar 1}{\bar 1}}$, and ${\omega}^{3}=
A_{11,1}$)

The reason for taking fractional exponents in (3.7) is to make our conditions
invariant under the scale change of contact form by a positive constant
multiple as we'll see later. Take a small amount of ${\int_M}
|E_{11,{\bar 1}1}|^{2}dv_{\theta}$ and ${\int_M}R^{2}|E_{11}|^{2} 
dv_{\theta}$ to gain the term ${\int_M}R|E_{11,{\bar 1}}|^{2}dv_{\theta}$
by (3.6) in the right-hand side of (3.5) while keeping the quadratic
hermitian form in $E_{11,{\bar 1}1}$, $iE_{11,0}$, $E_{11,1}$,
$E_{11,{\bar 1}}$, $E_{11}$ positive definite at least for $R=constant>0$,
$A_{11}=0$. For instance, we can take ${\lambda}={\rho}=\frac{1}{4}$ in
(3.6) and then use it and (3.7) in estimating the right-hand side of (3.5).
The final result reads 

\begin{eqnarray}
0&\geq&{\int_M}\{ \frac{29}{48}|E_{11,{\bar 1}1}|^{2}+2|E_{11,0}|^{2}
+\frac{1}{3}R|E_{11,1}|^{2}\\
 &  & +(\frac{1}{8}R-\frac{2}{3}|A_{11,1}|^{\frac{2}{3}})|E_{11,{\bar 1}}|^
{2}+[\frac{29}{48}R^{2}+\frac{11}{48}{\bigtriangleup}_{b}R
\nonumber \\
 &  & +6|A_{11}|^{2}+\frac{8i}{3}(A_{11,{\bar 1}{\bar 1}}-A_{{\bar 1}{\bar 1},
11})-\frac{2}{3}|A_{11,1}|^{\frac{4}{3}}]|E_{11}|^{2}
\nonumber \\
 & & +2Re[-iE_{11,0}E_{{\bar 1}{\bar 1},1{\bar 1}}+\frac{2}{3}iRE_{11,0}
E_{{\bar 1}{\bar 1}}-\frac{1}{6}RE_{11,{\bar 1}1}E_{{\bar 1}{\bar 1}}
\nonumber \\
 & & +(\frac{5}{48}R_{,{\bar 1}}-2iA_{{\bar 1}{\bar 1},1})E_{11,1}  
E_{{\bar 1}{\bar 1}}-\frac{5i}{3}A_{11}E_{11,{\bar 1}}E_{{\bar 1}{\bar 1},
{\bar 1}}]\} dv_{\theta}
\nonumber
\end{eqnarray}

The integrand in (3.8) is a quadratic hermitian form in $E_{11,{\bar 1}1}$,
$iE_{11,0}$, $E_{11,1}$, $E_{11,{\bar 1}}$, and $E_{11}$. By basic
linear algebra it is positive definite if and only if $R>0$,

\begin{equation}
{\left| \begin{array}{cccc}
\frac{29}{48} & -1 & 0 & 0 \\
-1 & 2 & 0 & 0 \\
0 & 0 & \frac{1}{3}R &\frac{5i}{3}A_{{\bar 1}{\bar 1}} \\
0 & 0 & \frac{-5i}{3}A_{11} & \frac{1}{8}R-\frac{2}{3}|A_{11,1}|^
{\frac{2}{3}}
\end{array}\right|}>0,\: and
\end{equation}

\begin{equation}
\left| \begin{array}{ccccc}
\frac{29}{48} & -1 & 0 & 0 & -\frac{1}{6}R \\
-1 & 2 & 0 & 0 & \frac{2}{3}R \\
0 & 0 & \frac{1}{3}R & \frac{5i}{3}A_{{\bar 1}{\bar 1}} &
\frac{5}{48}R_{,{\bar 1}}-2iA_{{\bar 1}{\bar 1},1} \\
0 & 0 & \frac{-5i}{3}A_{11} & \frac{1}{8}R-\frac{2}{3}|A_{11,1}|^
{\frac{2}{3}} & 0 \\
-\frac{1}{6}R & \frac{2}{3}R & \frac{5}{48}R_{,1}+2iA_{11,{\bar 1}} &
0 & \mathcal{M}
\end{array} \right|
\end{equation}

\noindent is larger than $0$.
Here $\mathcal{M}=\frac{29}{48}R^{2}+\frac{11}{48}{\bigtriangleup}
_{b}R+6|A_{11}|^{2}-\frac{2}{3}|A_{11,1}|^{\frac{4}{3}}+
\frac{8i}{3}(A_{11,{\bar 1}{\bar 1}}-A_{{\bar 1}{\bar 1},11})$.
A straightforward computation shows that (3.9) is equivalent to

\begin{equation}
\frac{3}{8}R^{2}-2R|A_{11,1}|^{\frac{2}{3}}-25|A_{11}|^{2}>0
\end{equation}

\noindent while (3.10) is equivalent to 

\begin{eqnarray}
&& (\frac{3}{8}R^{2}-2R|A_{11,1}|^{\frac{2}{3}}-25|A_{11}|^{2})\{
\frac{83}{3456}R^{2}\\
&& +\frac{55}{1152}{\bigtriangleup}_{b}R+\frac{5}{4}|A_{11}|^{2}
-\frac{5}{36}|A_{11,1}|^{\frac{4}{3}}+\frac{5i}{9}
(A_{11,{\bar 1}{\bar 1}}-A_{{\bar 1}{\bar 1},11})\}
\nonumber \\
&& -\frac{15}{8}(\frac{1}{8}R-\frac{2}{3}|A_{11,1}|^{\frac{2}{3}})
|\frac{5}{48}R_{,{\bar 1}}-2iA_{{\bar 1}{\bar 1},1}|^{2}>0
\nonumber
\end{eqnarray}

Observe that if $\theta$ changes by a positive constant multiple 
$k$, $R$ and $A_{11}$ change by multiplying $k^{-1}$ while 
$A_{11,1}$,$A_{{\bar 1}{\bar 1},1}$, and $R_{,{\bar 1}}$ change
by multiplying $k^{-\frac{3}{2}}$. Similarly $A_{11,{\bar 1}{\bar 1}}
$, $A_{{\bar 1}{\bar 1},11}$, and ${\bigtriangleup}_{b}R$ change
by multiplying $k^{-2}$. So the conditions (3.11),(3.12) are 
invariant under the change of contact form by a positive
constant multiple. Now $Theorem \: B$ follows from (3.8) under
the conditions (3.11),(3.12).

\end{document}